\documentclass[11pt,twoside, reqno]{article}

\usepackage{amssymb}
\usepackage{amsmath}
\usepackage{mathrsfs}
\usepackage{amsthm}
\usepackage{amsfonts}
\usepackage{txfonts}
\usepackage{enumerate}
\usepackage{hyperref}
\usepackage{graphicx}
\usepackage{caption2}
\usepackage{subfigure}
\usepackage{float}

\usepackage{latexsym,amssymb}
\usepackage{color}

\usepackage{graphicx}

\usepackage{pgf,tikz}
\usepackage{mathrsfs}
\usetikzlibrary{arrows}
\definecolor{ttttff}{rgb}{0.2,0.2,1.}
\definecolor{ttffcc}{rgb}{0.2,1.,0.8}
\definecolor{qqqqff}{rgb}{0.,0.,1.}

\definecolor{zzttqq}{rgb}{0.6,0.2,0.}
\definecolor{qqqqff}{rgb}{0.,0.,1.}

\allowdisplaybreaks

\usepackage{indentfirst}

\def\fz{\infty}

\def\r{\right}
\def\lf{\left}

\def\rr{{\mathbb R}}

\def\rn{{{\rr}^n}}

\def\cl{{\mathcal L}}

\def\dfrac{\displaystyle\frac}

\newtheorem{theorem}{Theorem}
\newtheorem{lemma}[theorem]{Lemma}

\newtheorem{proposition}[theorem]{Proposition}
\theoremstyle{definition}

\numberwithin{equation}{section}

\def\loc{{\mathop\mathrm{loc\,}}}

\numberwithin{equation}{section}

\begin{document}

\arraycolsep=1pt

\title{\bf\Large
Corrigendum to ``Hardy Spaces $H_\cl^1(\rn)$ Associated to
Schr\"odinger Type Operators $(-\Delta)^2+V^2$'' [Houston J. Math  36
(4) (2010),  1067-1095] \footnotetext{\hspace{-0.35cm} 2020 {\it
Mathematics Subject Classification}. Primary: 42B35; Secondary:
42B30, 42B25, 35J10.
\endgraf {\it Key words and phrases}.
Schr\"odinger type operator, heat kernel, Hardy spaces, atom, higher
order Riesz transform, fractional integral.
\endgraf
}}
\author{Jun Cao, Yu Liu and Dachun Yang\,\footnote{Corresponding author}}
\date{}
\maketitle

\vspace{-0.8cm}

\begin{center}
\begin{minipage}{13.5cm}
{\small {\bf Abstract.}  We rectify an incorrect citation of the reference in
obtaining the Gaussian upper bound for heat kernels of the Schr\"odinger
type operators $(-\Delta)^2+V^2$.}
\end{minipage}
\end{center}

Let $\cl:=(-\Delta^2)+V^2$ be a Schr\"odinger type operator on the Euclidean space
 $\rn$ with the potential $V\in L_\loc^1(\rn)$ for $n\ge 5$. Denote by
 $e^{-t\cl}(x,\,y)$ the heat kernel  of $\cl$  with $t\in(0,\,\fz)$. In our article \cite{cly10},
 we used in \cite[Lemma 2.4]{cly10} the following
\begin{lemma}\label{lem1}
 Let $V\in L_\loc^1(\rn)$. Then there exist positive constants $C$ and
 $A_5$ such that, for any $t\in (0,\fz)$ and $x$, $y\in \rn$,
 \begin{align}\label{1}
 \lf|e^{-t\cl}(x,y)\r|\le C
 t^{-n/4}\exp\lf\{-A_5\dfrac{|x-y|^{4/3}}{t^{1/3}}\r\}.
 \end{align}
 \end{lemma}
We claimed that the above lemma is  a simple corollary of  \cite[Proposition 5.2]{BaDa96},
where the latter says the following conclusion.
 \begin{proposition}
 Any heat kernel bound obtained for $(-\Delta)^2$ by means of Proposition
 2.5 therein is also valid for $\cl=(-\Delta)^2+V^2$.
 \end{proposition}

Although the heat kernel of the bi-Laplacian $(-\Delta)^2$ satisfies \eqref{1}
(see, for instance, \cite[Proposition 45]{AQ00}),  a dimension restriction
$n<4$ is implicitly used in  \cite[Proposition 5.2]{BaDa96}, while we worked
the case $n\ge 5$ in \cite{cly10}. Thus, Lemma \ref{lem1} is wrong.

To \emph{correct this mistake} which the aforementioned miscitation
brings, it is convenient to make the following assumption throughout
\cite{cly10}.

\bigskip

 \noindent{\textbf{Assumption (A).}}
 Let $V\in L_\loc^1(\rn)$. Assume that there exist positive constants $C$ and
 $A_5$ such that, for any $t\in (0,\fz)$ and $x$, $y\in \rn$,
 \begin{align*}
 \lf|e^{-t\cl}(x,y)\r|\le C t^{-n/4}\exp\lf\{-A_5\dfrac{|x-y|^{4/3}}{t^{1/3}}\r\}.
 \end{align*}

\medskip

We point out that this modification does not affect the main content 
of the article \cite{cly10}, because we worked mainly
on the real-variable theory of Hardy spaces associated with $\cl$, while the 
Gaussian upper bound serves only as the start point
of all arguments therein.  Note that, although Assumption (A) narrows down the scope of 
the potentials $V$ of the
Schr\"odinger type $(-\Delta)^2+V^2$ in \cite{cly10},
it is still meaningful, especially for those $V$  in some generalized Schechter-type
class (see \cite{clyz}).  For example, Assumption (A) holds true for any
$\cl=(-\Delta)^2+V^2$ with
\begin{align*}
V(x):=\lf[c(1+|x|)^a+1\right]^{1/2}, \ \forall\, x\in\rn,
\end{align*}
$a\in (-\fz,-4)$, and $c\in(0,\fz)$ being small enough. In this case,  $V^2$
belongs to the reverse H\"older class $B_\fz(\rn)$. Indeed, for any
$a\in (-\fz,-4)$, since $\|c(1+|x|)^a\|_{L^{n/4}(\rr^n)}$ is small
enough, it follows from \cite[Remark 5.9(ii)]{clyz} that the
Schr\"odinger type operator $(-\Delta)^2+c(1+|x|)^a$ satisfies
\eqref{1}. Then, using the Trotter product formula and noting that $1$ is
a positive constant, we conclude that
\begin{align*}
e^{-t[(-\Delta)^2+V^2]}=e^{-t[(-\Delta)^2+c(1+|x|)^a]}e^{-t},
\end{align*}
which implies that $(-\Delta)^2+V^2$ also satisfies \eqref{1}.

We point out that Assumption (A) may not hold true in the general case
$V\in L_\loc^1(\rn)$  and $n\ge 5$, as one can obtain only a local
version of the Gaussian upper bound for $V$ in an even smaller Kato
class (see, for instance, \cite{DeDiYa14}).

\bigskip

\noindent{\textbf{Acknowledgement.}} Dachun Yang would like to
thank Professor Jacek Dziubanski for pointing out this question to
him.

\bigskip

{\small\noindent Jun Cao

\noindent
Department of Applied Mathematics, Zhejiang University of Technology,
Hangzhou 310023, People's Republic of China

\smallskip

\noindent{\it E-mail:} \texttt{caojun1860@zjut.edu.cn}

\bigskip

\noindent {Yu Liu}
\smallskip

\noindent School of Mathematics and Physics, University of Science and Technology Beijing,
Beijing 100083, People's Republic of China

\smallskip

\noindent{\it E-mail:} \texttt{liuyu75@pku.org.cn}

\bigskip

\noindent {Dachun  Yang} (Corresponding author)
\smallskip

\noindent Laboratory of Mathematics and Complex Systems (Ministry of Education of China),
School of Mathematical Sciences, Beijing Normal University, Beijing 100875, 
People's Republic of China

\smallskip

\noindent{\it E-mail:} \texttt{dcyang@bnu.edu.cn}

\end{document}